\documentclass[12pt]{article}
\input epsf
\usepackage{amsmath}
\usepackage{amssymb}
\usepackage{theorem}

\numberwithin{equation}{section}

\newtheorem{thm}{Theorem}[section]
\newtheorem{lemma}[thm]{Lemma}
\newtheorem{prop}[thm]{Proposition}
\newtheorem{cor}[thm]{Corollary}
{\theorembodyfont{\rmfamily}
\newtheorem{defn}[thm]{Definition}

\newtheorem{problem}{Problem}

\newtheorem{rmk}[thm]{Remark}
}

\newcommand{\qed}{\hfill \mbox{\raggedright \rule{.07in}{.1in}}}
 
\newenvironment{proof}{\vspace{1ex}\noindent{\bf
Proof}\hspace{0.5em}}{\hfill\qed\vspace{1ex}}
\newenvironment{pfof}[1]{\vspace{1ex}\noindent{\bf Proof of
#1}\hspace{0.5em}}{\hfill\qed\vspace{1ex}}

\newcommand{\R}{{\mathbb R}}
\newcommand{\C}{{\mathbb C}}
\newcommand{\Z}{{\mathbb Z}}
 \newcommand{\N}{{\mathbb N}}
 \newcommand{\T}{{\mathbb T}}

\newcommand{\diam}{\operatorname{diam}}

\newcommand{\SMALL}{\textstyle}

\title{A Vector-Valued Almost Sure Invariance Principle \\ for
Hyperbolic Dynamical Systems}

\author{
Ian Melbourne \thanks{Department of Mathematics and Statistics, University of Surrey,
Guildford GU2 7XH, UK.  E-mail: ism@math.uh.edu}
\and 
Matthew Nicol \thanks{ Department of Mathematics, University of Houston,
Houston TX 77204-3008, USA.  Email: nicol@math.uh.edu}
}

\date{20 June, 2006}

\begin{document}

\maketitle

\begin{abstract}
We prove an almost sure invariance principle (approximation by $d$-dimensional
Brownian motion) for vector-valued H\"older observables of large classes of
nonuniformly hyperbolic dynamical systems.   
These systems include Axiom~A diffeomorphisms and flows as well as 
systems modelled by Young towers with moderate tail decay rates.

In particular, the position variable of the planar periodic Lorentz gas
with finite horizon approximates a $2$-dimensional Brownian motion.
\end{abstract}

\section{Introduction}

The scalar almost sure invariance principle (ASIP), or approximation
by one-dimensional Brownian motion, is a strong
statistical property of sequences 
of random variables introduced by Strassen~\cite{Strassen64,Strassen67}.
It implies numerous other statistical limit laws
including the central limit theorem, the functional central limit theorem,
and the law of the iterated logarithm.  See~\cite{HallHeyde80,PhilippStout75}
and references therein for a survey of consequences of the ASIP.

The scalar ASIP has been shown to hold for large classes of dynamical 
systems~\cite{ConzeBorgne01,DenkerPhilipp84,Dolgopyat04,FMT03,HofbauerKeller82,HMsub,MN05,MT02,Nagayama04}.
Chernov \& Dolgopyat~\cite[Problem~1]{ChernovDolgopyatapp} asked for a proof of
the ASIP for $\R^d$-valued observables, and it is this problem that is solved 
in this paper.
Our main result applies to a large variety of dynamical systems, as surveyed 
in Section~\ref{sec-app}.  

As a secondary matter, we obtain explicit error estimates that depend on
the dimension $d$ and the lack of hyperbolicity.   Even for $d=1$, this 
estimate is better than those in almost all of the above references.
The exception is~\cite{FMT03} which gives the best available estimate for
scalar ASIPs for a restricted class of systems.

\subsection{Statement of the main results}

\begin{defn} \label{def-ASIP}
A sequence $\{S_N\}$ of random variables with values in $\R^d$
satisfies a {\em $d$-dimensional almost sure invariance 
principle (ASIP)} if there exists $\lambda>0$ and a probability space 
supporting a sequence of random variables 
$S_N^*$ and a $d$-dimensional Brownian motion $W(t)$ such that
\begin{itemize}
\item[(a)] $\{S_N;N\ge1\}=_d\{S_N^*;N\ge1\}$, and
\item[(b)] $S_N^*=W(N)+O(N^{\frac12-\lambda})$ as $N\to\infty$ almost 
everywhere.
\end{itemize}
For brevity, we write
$S_N=W(N)+O(N^{\frac12-\lambda})\enspace\text{a.e.}$
The ASIP for a one-parameter family $S_T$ of $\R^d$-valued random variables
is defined similarly, and denoted 
$S_T=W(T)+O(T^{\frac12-\lambda})\enspace\text{a.e.}$
\end{defn}

\begin{rmk}   The ASIP is said to be {\em nondegenerate} if the Brownian motion
$W(t)$ has nonsingular covariance matrix $\Sigma$.
For the classes of dynamical systems considered in this paper, the
ASIP is nondegenerate for {\em typical} observables.  More precisely, there
is a closed subspace $Z$ of infinite codimension in the space of
all (piecewise) H\"older $\R^d$-valued observables such that $\Sigma$
is nonsingular whenever $\phi\not\in Z$.  (By considering all
one-dimensional projections it suffices to consider the case $d=1$.
This is done explicitly in for example~\cite[Section~4.3]{HMsub}.)
\end{rmk}

\paragraph{Axiom~A diffeomorphisms and flows}
Our results are most easily stated in the uniformly hyperbolic (Axiom~A) 
context.

\begin{thm}  \label{thm-A}
Let $f:M\to M$ be a diffeomorphism
with a (nontrivial) uniformly hyperbolic basic set $X\subset M$,
and suppose that $\mu$ is an equilibrium measure corresponding to a H\"older 
potential.  Let $\phi:X\to\R^d$ be a mean zero H\"older observable with 
partial sums $S_N=\sum_{n=1}^N\phi\circ f^j$.  Then for any $\epsilon>0$,
\[
S_N=W(N)+O(N^{\beta+\epsilon})\enspace\text{a.e.}
\quad\text{where}\enspace \beta=\SMALL\frac{2d+3}{4d+7}. 
\]
(The improved estimate $\beta=\frac14$ holds when $d=1$~\cite{FMT03}.)
\end{thm}

An immediate consequence (see for example~\cite{DenkerPhilipp84,MT04})
is the corresponding result for Axiom~A flows.

\begin{cor}  \label{cor-A}
Let $f_t:M\to M$ be a smooth flow with a (nontrivial) uniformly hyperbolic 
basic set $X\subset M$, and suppose that $\mu$ is an equilibrium measure 
corresponding to a H\"older potential.  
Let $\phi:X\to\R^d$ be a mean zero H\"older observable with partial sums 
$S_T=\int_0^T\phi\circ f_t\,dt$.  Then for any $\epsilon>0$,
\[
S_T=W(T)+O(T^{\beta+\epsilon})\enspace\text{a.e.}
\quad\text{where}\enspace \beta=\SMALL\frac{2d+3}{4d+7}. 
\]
(The improved estimate $\beta=\frac14$ holds when $d=1$~\cite{FMT03,MT04}.)
\qed
\end{cor}

\begin{rmk}
Denker \& Philipp~\cite{DenkerPhilipp84} proved Theorem~\ref{thm-A}
and Corollary~\ref{cor-A} in the case $d=1$ (though with a weaker
error term).  
\end{rmk}

\paragraph{Nonuniformly hyperbolic systems}
Our results apply also to maps $f:M\to M$ that
are nonuniformly expanding/hyperbolic in the sense of 
Young~\cite{Young98,Young99}.   Roughly speaking, such maps possess a
subset $\Lambda\subset M$ 
and a return time $R:\Lambda\to\Z^+$
such that the induced map $f^R:\Lambda\to \Lambda$ is uniformly hyperbolic.

\begin{thm} \label{thm-NUH}   Let $f:M\to M$ be a diffeomorphism
(possibly with singularities) that is nonuniformly hyperbolic
in the sense of Young~\cite{Young98,Young99}.
In particular, $f$ satisfies conditions~(A1)--(A4) in
Section~\ref{sec-NUH} and possesses an SRB measure $m$.
Assume that the return time function $R$ lies in $L^p$, $p>2$.
Let $\phi:M\to\R^d$ be a mean zero H\"older observation
with partial sums $S_N=\sum_{n=1}^N\phi\circ f^j$.
Then for any $\epsilon>0$,
\[
S_N=W(N)+O(N^{\beta+\epsilon})\enspace\text{a.e.}
\quad\text{where}\enspace \beta=\SMALL\frac{\frac1p+2d+3}{4d+7}.
\]
If $d=1$ then $\beta$ can be improved to
$\beta=\frac{1}{2p}+\frac14$ for $2<p\le 4$ and
$\beta=\frac38$ for $p\ge4$.
\end{thm}

Again, there is an immediate extension to nonuniformly hyperbolic flows.
Suppose that $f:M\to M$ satisfies the assumptions of Theorem~\ref{thm-NUH}
with $R\in L^p$, $p>2$, and that $f_t$ is a suspension flow over $f$ with
a (uniformly bounded) H\"older roof function.   By~\cite{MT04},
$\R^d$-valued H\"older observables of the suspension flow satisfy an ASIP of 
the form $S_T=W(T)+O(T^{\beta+\epsilon})$ a.e. where $\beta$ is as
in Theorem~\ref{thm-NUH}.

\paragraph{Application to Lorentz gases}
The planar periodic Lorentz gas was introduced by Sinai~\cite{Sinai70}.
This is a three-dimensional flow with phase space
$(\R^2-\Omega)\times S^1$, where $\Omega\subset\R^2$ is a 
periodic array of disjoint convex regions with $C^3$ boundaries.  
The coordinates are position $q\in \R^2-\Omega$ and velocity $v\in S^1$.
The flow satisfies
the finite horizon condition if the time between collisions with
$\partial\Omega$ is uniformly bounded.   

Let $q(t)\in\R^2$ denote the position at time $t$ of a particle starting
at position $q(0)$ pointing in direction $v(0)$.
Bunimovich \& Sinai~\cite{BunimSinai81}, see also~\cite{BunimSinaiChernov91},
proved that $q(t)$ satisfies a two-dimensional
functional central limit theorem (weak invariance principle)
supporting the view of such flows as a
deterministic model for Brownian motion.  
We complete this circle of ideas by proving the strong version of this result.

\begin{thm}\label{thm-lorentz}
Consider a planar periodic Lorentz gas satisfying the finite horizon condition.
Let $\epsilon>0$.  There is a two-dimensional Brownian motion $W(t)$
with nonsingular covariance matrix such that for almost every initial 
condition, $q(T)=W(T)+O(T^{\frac{7}{15}+\epsilon})$.
\end{thm}

\begin{rmk}
(a)  A number of authors~\cite{Chernovapp,MN05,Nagayama04} have independently 
established scalar ASIPs for one-dimensional projections of $q(t)$.
In hindsight, the scalar ASIP for the Lorentz gas follows from
earlier work of~\cite{FMT03}, 
again with $\beta=\frac14$.
(We note that the methods of~\cite{FMT03} apply in the first place only to 
the time-reversal of the dynamical system under study.  Their
applicability here is due to the fact that the class of systems is closed 
under time-reversal.)
\\[.75ex]
(b) The finite horizon condition is crucial.   For infinite horizons,
Sz\'asz \& Varj\'u~\cite{SzaszVarjusub}
prove that $q(t)$ lies in the nonstandard domain of the normal distribution.
In particular, the central limit theorem fails, hence the ASIP fails.
\end{rmk}

\subsection{Consequences of the vector-valued ASIP}

For convenience, we suppose that the Brownian motion in the ASIP
is nondegenerate.  Coordinates can be chosen on $\R^d$ so that
$W(t)$ is a standard $d$-dimensional Brownian motion with $\Sigma=I_d$.
Throughout, the norm on $\R^d$ is taken to be the usual Euclidean norm.
The following consequences of the ASIP are summarised 
in~\cite[p.~233]{MorrowPhilipp82}.
Here, LIL stands for {\em law of the iterated logarithm} and the functional
LIL stated below is a far-reaching generalisation, due to Strassen, of
the classical LIL.   

\begin{prop} For the dynamical systems to which the results in this paper apply,
the following consequences hold (after normalisation so that $\Sigma=I_d$):

\vspace{1ex}
\noindent$\bullet$ {\bf Functional LIL}
Let $C=C([0,1],\R^d)$ be the Banach space of continuous maps $f:[0,1]\to\R^d$
with the supremum norm.    Let $K$ be the (compact) set of $f\in C$ absolutely
continuous with $f(0)=0$, $\int_0^1|f'(t)|^2dt\le1$.   
Define $f_n(i/n)=S_i/\sqrt{2n\log\log n}$, $i=0,\dots,n$, and linearly 
interpolate to obtain $f_n\in C$.  Then a.s.\ the sequence
$\{f_n\}$ is relatively compact in $C$ and its set of limit points is
precisely $K$.

\vspace{1ex}
\noindent$\bullet$ {\bf Upper and lower class refinement of the LIL}
Let $\phi(t):\R\to\R$ be a positive
nondecreasing function.   Then 
\[
P(|S_N|>N^{\frac12}\phi(N)\enspace\text{i.o.})=0 \enspace\text{or}\enspace 1
\]
according to whether
$\int_1^\infty \frac{\phi^d(u)}{u}\exp(-\frac12\phi^2(u))\,du$ converges or diverges.

\vspace{1ex}
\noindent$\bullet$ {\bf Upper and lower class refinement of Chung's LIL}
Let $\phi(t):\R\to\R$ be a positive
nondecreasing function.   Then there is a constant $c$ (depending only on $d$)
such that
\[
P(\max_{n\le N}|S_n|<cN^{\frac12}\phi^{-1}(N)\enspace\text{i.o.})=0 \enspace\text{or}\enspace 1
\]
according to whether
$\int_1^\infty \frac{\phi^2(u)}{u}\exp(-\phi^2(u))\,du$ converges or diverges.

\vspace{1ex}
\noindent$\bullet$ {\bf Central limit theorem and functional central limit theorem}
\qed
\end{prop}

\begin{rmk} \label{rmk-Berger}
 (a) Berger~\cite{Berger90} gives a unified approach to the ASIP
for weakly dependent sequences of random variables with values in a
 real separable Banach space, but with error term $o(\sqrt{N\log\log N})$.
It follows from Berger~\cite[Corollary~4.1, part A.5]{Berger90} 
and Melbourne \& Nicol~\cite{MN05} that the Banach space-valued ASIP formulated
in~\cite[Theorem~3.2]{Berger90} holds for all dynamical systems considered in 
this paper.   In particular, the $\R^d$-valued ASIP holds with error term 
$o(\sqrt{N\log\log N})$.
This error term suffices for the functional LIL, but is inadequate for the 
upper and lower class refinements and for the (functional) central limit 
theorem; whereas the error term established in this paper suffices.
Indeed this was the original motivation of 
Jain {\em et al.}~\cite{JainJogdeoStout75}
to improve the error term in Strassen's scalar ASIP.
\\[.75ex]
(b)  The $\R^d$-valued functional central limit theorem, being a 
distributional result, can be proved directly under the more general 
condition $R\in L^2$ in Theorem~\ref{thm-NUH}: reduce as in this paper to the 
setting in Section~\ref{sec-GM} and then apply the 
method of~\cite[Section~3.3]{FMT03}.
\end{rmk}

We end this section by discussing briefly the probabilistic
methods used in this paper.
Strassen's original proof of the scalar ASIP for IIDs and martingales~\cite{Strassen64,Strassen67} relies heavily on the Skorokhod embedding theorem
for scalar stochastic processes.   
This method was extended to {\em weakly dependent} sequences of random variables
by a number of authors, using {\em blocking arguments} to
reduce to the martingale case, see~\cite{PhilippStout75}.
In particular, Philipp \& Stout~\cite[Theorem~7.1]{PhilippStout75} formulated
a version of the scalar ASIP which is particularly useful for dynamical
systems~\cite{HofbauerKeller82,DenkerPhilipp84,MN05}.  

Attempts to extend Strassen's proof to the $\R^d$-valued situation were 
only partially successful~\cite{Kiefer72}, and the same is true
for the completely different {\em quantile transform method}
of Cs\"org\"o \& R\'ev\'esz~\cite{CsorgoRevesz75}.
Eventually, Berkes \& Philipp~\cite{BerkesPhilipp79} introduced a third
method which works in any number of dimensions, and the
applicability of this method was extended to weakly dependent sequences by
Kuelbs \& Philipp~\cite{KuelbsPhilipp80}.

The remainder of the paper is organised as follows.
In Section~\ref{sec-ASIP}, we combine the blocking argument 
in~\cite{PhilippStout75} 
with the results of~\cite{BerkesPhilipp79,KuelbsPhilipp80} to
prove a vector-valued ASIP for $\R^d$-valued 
random variables satisfying certain hypotheses.   
In Section~\ref{sec-GM}, we verify these hypotheses for Gibbs-Markov maps
and derive Theorem~\ref{thm-A} as a consequence.
In Section~\ref{sec-app}, we first prove the ASIP for nonuniformly expanding
maps and then prove Theorems~\ref{thm-NUH} and~\ref{thm-lorentz}.
We also list numerous other situations to which
our results apply, and we mention some open problems regarding time-one
maps of flows.

\section{A vector-valued ASIP for functions of mixing sequences}
\label{sec-ASIP}

In this section, we prove a vector-valued ASIP for $\R^d$-valued 
random variables satisfying certain hypotheses.   
This is the vector-valued analogue of~\cite[Theorem~7.1]{PhilippStout75}
though with hypotheses tailored to the dynamical systems setting.
(A result of this type is hinted at in Kuelbs \& 
Philipp~\cite{KuelbsPhilipp80}, but it is necessary to work through the details
to determine the hypotheses, which
were left unstated.  In any case, the estimates  in~\eqref{eq-Burk}
and~\eqref{eq-CLT} are not so natural in the probabilistic setting 
in~\cite{KuelbsPhilipp80}, and partly account for our strong error term.)

\subsection{Statement of the ASIP}

Let $\xi_1,\xi_2,\ldots$ be a 
sequence of real-valued random variables 
and let $\mathcal{F}_a^b=\sigma\{\xi_n;a\le n\le b\}$.
We assume the strong-mixing condition
\begin{align} \label{eq-mix}
|P(AB)-P(A)P(B)|\le C\tau^n \quad\text{for all $A\in\mathcal{F}_1^k$
and $B\in\mathcal{F}_{k+n}^\infty$.}
\end{align}

Let $p>2$, and let $\eta_n$ be a strictly stationary sequence of 
$\mathcal{F}_n^\infty$-measurable $\R^d$-valued random variables 
satisfying
\begin{align} \label{eq-Lp}
E\eta_n=0\quad\text{and}\quad |\eta_n|_p=C,
\end{align}
and the (backwards) Burkholder-type inequality
\begin{align}  \label{eq-Burk}
\Bigl|\max_{1\le \ell\le N}\,\bigl|{\SMALL\sum_{n=\ell}^N}\,\eta_n\bigr|\,\Bigr|_p\le 
C N^{\frac12}. 
\end{align}
Define $\eta_{\ell n}=E(\eta_n|\mathcal{F}_n^{n+\ell})$.
We require that
\begin{align} \label{eq-local}
|\eta_n-\eta_{\ell n}|_p\le C\tau^\ell.
\end{align}

Let $\Sigma$ be a symmetric positive semidefinite $d\times d$ covariance
matrix.
Given $u\in\R^d$, define
$f_N(u)=E\exp(i\langle u, \sum_{n\le N}\eta_n/\sqrt N\rangle)$ and 
$g(u)=\exp(-\frac12\langle u,\Sigma u\rangle)$.
Assume that there exists $\epsilon>0$ such that
\begin{align} \label{eq-CLT}
|f_N(u)-g(u)|\le CN^{-\frac12} 
\enspace\text{for all $|u|\le \epsilon N^{\frac12}$}.
\end{align}

\begin{thm} \label{thm-ASIP}
Assume conditions~\eqref{eq-mix}--\eqref{eq-CLT}.
Let 
$\beta=\frac{\frac1p+2d+3}{4d+7}\in\Bigl[\frac{2d+3}{4d+7},\frac12\Bigr)$,
and let $\epsilon>0$.
Then 
there is a $d$-dimensional Brownian motion $W(t)$ with covariance
matrix $\Sigma$ such that 
$\sum_{n\le N}\eta_n = W(N)+O(N^{\beta+\epsilon})\enspace\text{a.e.}$
\end{thm}

\begin{rmk}   (a) For $d=1$, we obtain under similar hypotheses,
but using a different method, the improved error estimate 
$\beta=\frac{1}{2p}+\frac14$ for $2<p\le 4$ and
$\beta=\frac38$ for $p\ge4$.
See Appendix~\ref{sec-1}.
\\[.75ex]
(b)  It is evident from the proof of Theorem~\ref{thm-ASIP} that
the exponential rates in~\eqref{eq-mix} and~\eqref{eq-local}
can be replaced by sufficiently high polynomial rates.
Further relaxing of the assumptions is possible at the cost of obtaining
a weaker estimate in Theorem~\ref{thm-ASIP}.
\end{rmk}

\subsection{Preliminaries}

\indent The following result of~\cite{Davydov70,VolkonskiiRozanov59}
is stated as~\cite[Lemma~7.2.1]{PhilippStout75}

\begin{lemma} \label{lem-Dav}
Let $\mathcal{F}, \mathcal{G}$ be $\sigma$-fields and $\beta\ge0$ such that
$|P(AB)-P(A)P(B)|\le \beta$ for all $A\in\mathcal{F}$, $B\in\mathcal{G}$.
Let $p,q,r>1$ satisfy $\frac1p+\frac1r+\frac1s=1$.  
Suppose that $\xi\in L^r(\mathcal{F})$, $\eta\in L^s(\mathcal{G})$.
Then
$|E(\xi\eta)-E(\xi)E(\eta)|\le 10\beta^\frac1p\|\xi\|_r\|\eta\|_s$.
\end{lemma}

The next result is due to Dvoretsky~\cite{Dvoretzky72}, 
see~\cite[Lemma~2.2]{KuelbsPhilipp80}.

\begin{lemma} \label{lem-Dv}
Let $\mathcal{F},\mathcal{G}$ be $\sigma$-fields.
Let $\xi$ be a complex-valued $\mathcal{F}$-measurable
random variable with $|\xi|\le1$.
Then
$E|E(\xi|\mathcal{G})-E\xi|\le 2\pi\sup_{A\in\mathcal{F},B\in\mathcal{G}}
|P(AB)-P(A)P(B)|$.
\end{lemma}

The following Gal-Koksma strong law~\cite{GalKoksma50} 
is stated in~\cite[Theorem~A1]{PhilippStout75}.

\begin{lemma} \label{lem-GK}
Let $\xi_j$ be a sequence of random variables with $E\xi_j=0$, and let $q>0$.
Suppose that $E|\sum_{j=m}^n \xi_j|^2\le n^q-m^q$
for all $n\ge m\ge1$.   For any
$\epsilon>0$, $\sum_{j=1}^M\xi_j\ll M^{\frac{q}{2}+\epsilon}$ a.e.
\end{lemma}

\subsection{Introduction of the blocks}
\label{sec-block}

Fix $Q>\alpha>0$.   
Define random variables $y_1,z_1,y_2,z_2\ldots$
consisting of sums of consecutive $\eta_{\ell(n),n}$ where the $j$'th blocks
$y_j$ and $z_j$ consist of $[j^Q]$ and $[j^{\alpha}]$ such terms
respectively, and throughout the $j$'th blocks
$\ell(n)=[\frac12 j^\alpha]$.  

In other words, $y_j=\sum_n\eta_{\ell,n}$,
where $\ell=[\frac12 j^\alpha]$ and the sum ranges over 
$\sum_{i=1}^{j-1}([i^Q]+[i^\alpha])<n\le 
\sum_{i=1}^{j-1}([i^Q]+[i^\alpha])+[j^Q]$.
Similarly for $z_j$.

Let $\mathcal{L}_a^b=\sigma\{y_j;a\le j\le b\}$
and $\tilde{\mathcal{L}}_a^b=\sigma\{z_j;a\le j\le b\}$.

\begin{lemma} \label{lem-mix}
There exists (a modified) $\tau\in(0,1)$ such that for all 
$k,n\ge1$,
\[
|P(AB)-P(A)P(B)|\ll \tau^{(k+n)^\alpha} \quad
\text{for all $A\in\mathcal{L}_1^k$ and $B\in\mathcal{L}_{k+n}^\infty$.}
\]
The same is true for all $A\in\tilde{\mathcal{L}}_1^k$ and $B\in\tilde{\mathcal{L}}_{k+n}^\infty$.
\end{lemma}

\begin{proof}
Note that $\mathcal{L}_1^k$ is defined using $y_1,\ldots,y_k$
which are defined using $\eta_{\ell n}$ with
$\ell\le[\frac12 k^{\alpha}]$, $n\le \sum_{i=1}^{k-1}([i^Q]+[i^{\alpha}])+[k^Q]$.    This involves 
conditioning on $\xi_n$ with $n\le \sum_{i=1}^{k-1}([i^Q]+[i^{\alpha}])+[k^Q]+[\frac12 k^{\alpha}]$.
Similarly for $\mathcal{L}_{k+n}^\infty$ and we obtain
\[
\mathcal{L}_1^k\subset\mathcal{F}_1^{\sum_{i=1}^{k-1}([i^Q]+[i^{\alpha}])+[k^Q]+[\frac12 k^{\alpha}]},
\quad
\mathcal{L}_{k+n}^\infty\subset\mathcal{F}_{\sum_{i=1}^{k+n-1}([i^Q]+[i^{\alpha}])+1}^\infty.
\]
Hence $|P(AB)-P(A)P(B)|\le \tau^N$ where
$N=\sum_{i=k+1}^{k+n-1}([i^Q]+[i^{\alpha}])+[k^{\alpha}]
-[{\SMALL\frac12} k^{\alpha}]+1$.
For all $k,n\ge1$, we compute that $N\gg (k+n)^\alpha$ as required for the
first statement.  (Note that the details for the cases $n=1$ and $n\ge2$ are 
slightly different.)
The second statement is proved in the same way.
\end{proof}

For $N\ge1$, let $y_{M_N},z_{M_N}$ be the pair of blocks that contains 
$\eta_{\ell(N),N}$.
Write 
\[
\SMALL y_{M_N}+z_{M_N}=\sum_{j=P_{M_N-1}+1}^{P_{M_N}} \eta_{\ell j}, \quad 
\ell=[{\SMALL\frac12} M_N^{\alpha}].
\]
In particular, $P_{M_N-1}<N\le P_{M_N}$, and 
$P_{M_N}-P_{M_N-1}=[M_N^Q]+[M_N^\alpha]\sim M_N^Q$.
It is immediate that

\begin{prop} \label{prop-M}
Writing $M=M_N$, we have $M^{1+Q} \sim {\SMALL\sum_{j\le M}}j^Q\sim N$.
In particular, $P_M-P_{M-1}\sim N^{Q/(1+Q)}$.\qed
\end{prop}

\begin{prop} \label{prop-local}
$\sum_{n\ge1}|\eta_n-\eta_{\ell n}|_p<\infty$.
\end{prop}

\begin{proof}
Focusing on the $M$'th block, and applying~\eqref{eq-local},
we obtain
$\sum_{P_{M-1}<n\le P_M}|\eta_n-\eta_{\ell(n),n}|_p\ll
M^Q\,\tau^{\frac12 M^\alpha}$ which is summable.
\end{proof}

\begin{prop} \label{prop-yz}
$|y_j|_p\ll j^{\frac12 Q}$ and $|z_j|_p\ll j^{\frac12\alpha}$.
\end{prop}

\begin{proof}
Write $y_j=\sum^*\eta_{\ell n}$ 
where $\sum^*=\sum_{n=a_j+1}^{a_j+[j^Q]}$.
By Proposition~\ref{prop-local},~\eqref{eq-Burk} and stationarity,
$|y_j|_p\le |\sum^*(\eta_{\ell n}-\eta_n)|_p
+ |\sum^*\eta_n|_p
\ll 1+ j^{\frac12 Q}\ll j^{\frac12 Q}$.
Similarly for $z_j$.~
\end{proof}

\subsection{Approximation result}

\begin{thm} \label{thm-approx}
Let $\beta=\max\{\frac1p+\frac12 Q,\frac12\}/(1+Q)$.
For any $\epsilon>0$, there exists $\alpha>0$ such that 
$\sum_{n\le N}\eta_n-\sum_{j\le M_N}y_j\ll N^{\beta+\epsilon}$ a.e.
\end{thm}

Begin by writing
\[
\sum_{n\le N}\eta_n -\sum_{j\le M_N}y_j=
\Bigl(\sum_{n\le P_{M_N}}\eta_n -\sum_{j\le M_N}(y_j+z_j)\Bigr)
 -\sum_{n=N+1}^{P_{M_N}}\eta_n +\sum_{j\le M_N}z_j.
\]
In the next three lemmas, we estimate these three terms
(following~\cite[Lemmas~7.3.2, 7.3.3, 7.3.4]{PhilippStout75}).
The result follows by combining these estimates.

\begin{lemma} \label{lem-blocks1}
$\sum_{n\le P_{M_N}}\eta_n -\sum_{j\le M_N}(y_j+z_j) \ll 1 \enspace\text{a.e.}$
\end{lemma}

\begin{proof}
By Proposition~\ref{prop-local}, 
$\sum_{n\le \infty}|\eta_n-\eta_{\ell n}|<\infty$ a.e.
Hence $|\sum_{n\le P_M}\eta_n-\sum_{j\le M}(y_j+z_j)|
=|\sum_{n\le P_M}(\eta_n-\eta_{\ell n})|\le 
\sum_{n\le \infty}|\eta_n-\eta_{\ell n}|\ll1$ a.e.
\end{proof}

\begin{lemma} \label{lem-blocks2}
Let $\beta=(\frac12+\frac12\alpha)/(1+Q)$.  For any $\epsilon>0$,
$\sum_{j\le M_N}z_j \ll N^{\beta+\epsilon} \enspace\text{a.e.}$
\end{lemma}

\begin{proof}
We have $|z_j|_{p}\ll j^{\frac12\alpha}$ and so 
$\sum_m^nEz_j^2\ll \sum_m^n j^{\alpha}\ll 
n^{1+\alpha}-m^{1+\alpha}$.
By Lemmas~\ref{lem-Dav} and~\ref{lem-mix}
(with $\tilde\tau=\tau^{\epsilon}$ where $\epsilon=1-2/p$), for all $i<j$,
\begin{align*}
|E z_iz_j| & \ll 
|z_i|_p|z_j|_p\tilde\tau^{j^\alpha} 
\le (i^\alpha\tilde\tau^{i^\alpha}j^\alpha\tilde\tau^{j^\alpha})^{\frac12}
\end{align*}
which is summable over $(i,j)\in\N^2$.  We have shown that 
$E(\sum_{j=m}^nz_j)^2 \ll n^{1+\alpha}-m^{1+\alpha}$,
for all $1\le m\le n$.  By Lemma~\ref{lem-GK},
$\sum_{j\le M}z_j\ll M^{\gamma}$ a.e. for $\gamma>\frac{1}{2}(1+\alpha)$,
and the result follows from Proposition~\ref{prop-M}.
\end{proof}

\begin{lemma} \label{lem-blocks3}
Let $\beta=(\frac1p+\frac12 Q)/(1+Q)$.   For any $\epsilon>0$,
$\sum_{n=N+1}^{P_{M_N}}\eta_n \ll N^{\beta+\epsilon} \enspace\text{a.e.}$
\end{lemma}

\begin{proof}
Let $A_M=\max_{P_{M-1}+1\le N\le P_M}|\sum_{n=N+1}^{P_M}\eta_n|$.
By~\eqref{eq-Burk} and stationarity, 
$|A_M|_p\ll (P_M-P_{M-1})^\frac12\ll M^{\frac12 Q}$.  Hence
\[
P(A_M>M^{\gamma})= P(A_M^p>M^{p\gamma}) \ll M^{-p(\gamma-\frac12 Q)},
\]
which is summable provided $\gamma>\frac1p+\frac12 Q$.   By Borel-Cantelli,
$A_M \ll M^\gamma$ a.e. and the result follows
from Proposition~\ref{prop-M}.
\end{proof}


\subsection{Proof of Theorem~\ref{thm-ASIP}}

We follow the argument of Kuelbs \& Philipp~\cite{KuelbsPhilipp80}
which extends Berkes \& Philipp~\cite{BerkesPhilipp79}.
Let $X_j=[j^Q]^{-\frac12}y_j$.   Note that $\mathcal{L}_1^j$ is an
increasing sequence of $\sigma$-fields such that $X_j$ is 
$\mathcal{L}_1^j$-measurable.

\begin{prop} \label{prop-BP1}
Let $\gamma\in(0,\frac12 Q)$.   There exists $\epsilon>0$ such that
$E\bigl|E(\exp(i\langle u,X_j\rangle)|\mathcal{L}_1^{j-1})-\exp(-\frac12\langle u,\Sigma u\rangle)\bigr|\le C' j^{\gamma-\frac12 Q}$
for all $u\in\R^d$ satisfying $|u|\le\epsilon j^\gamma$.
\end{prop}

\begin{proof}
Let $f_N(u)=E\exp(i\langle u, \sum_{n\le N}\eta_n/\sqrt N\rangle)$,
$g(u)=\exp(-\frac12\langle u,\Sigma u\rangle)$, 
and write 
\begin{align*}
& E\{\exp(i\langle u,X_j\rangle)|\mathcal{L}_1^{j-1}\}-g(u) = 
\bigl(E\{\exp(i\langle u,X_j\rangle)|\mathcal{L}_1^{j-1}\}- E\exp(i\langle u,X_j\rangle)\bigr) \\ & 
+\Bigl(E\exp(i\langle u,[j^Q]^{-\frac12}y_j\rangle)-
E\exp(i\langle u,[j^Q]^{-\frac12}{\SMALL\sum_{n\le [j^Q]}}\eta_n\rangle)\bigr) 
+\bigl(f_{[j^Q]}(u)-g(u)\bigr) \\ & = I+II+III.
\end{align*}

Using Lemmas~\ref{lem-Dv} and~\ref{lem-mix},
$E|I|\ll \tau^{j^\alpha}$.
Also, $III$ is estimated by~\eqref{eq-CLT} so it remains to estimate $II$.
Write $y_j=\sum^*\eta_{\ell n}$ where $\sum^*=\sum_{n=a_j+1}^{a_j+[j^Q]}$.
By stationarity and Proposition~\ref{prop-local},
\begin{align*}
|II| & =|E\Bigl(\exp(i\langle u,[j^Q]^{-\frac12}{\SMALL\sum^*}\eta_{\ell n}\rangle)-
\exp(i\langle u,[j^Q]^{-\frac12}{\SMALL\sum^*}\eta_n\rangle)\bigr)| \\
& \le|\exp(i\langle u,[j^Q]^{-\frac12}{\SMALL\sum^*}(\eta_{\ell n}-\eta_n)\rangle) 
-1|_1 \le
|\langle u,[j^Q]^{-\frac12}{\SMALL\sum^*}(\eta_{\ell n}-\eta_n)\rangle) |_1
\\ & \SMALL \le \epsilon j^\gamma [j^Q]^{-\frac12}|\sum_{n\ge1}(\eta_{\ell n}-\eta_n)|_1
\ll j^{\gamma-\frac12 Q}
\end{align*}
as required.
\end{proof}

\begin{prop} \label{prop-BP2}
Let $G$ be the distribution function of $N(0,\Sigma)$.
Then $G\{u:|u|>T\}\le e^{-DT^2}$.
\end{prop}

\begin{proof}  This is a straightforward calculation,
see for example~\cite[p.~43]{BerkesPhilipp79}.
\end{proof}

Let $\lambda_j=C'j^{\gamma-\frac12 Q}$, $T_j=\epsilon j^\gamma$,
where $\gamma\in(0,\frac12 Q)$ is chosen below.
By Propositions~\ref{prop-BP1} and~\ref{prop-BP2}, we have
\begin{align*}
& E|E\{\exp(i\langle u,X_j\rangle)|\mathcal{L}_1^{j-1}\}-g(u)|\le \lambda_j
\enspace\text{for all
$|u|\le T_j$}, \\[.75ex]
& G\{u:|u|>{\SMALL\frac14} T_j\}\le \delta_j,
\end{align*}
where $\delta_j=e^{-D'j^{2\gamma}}$.
These are the hypotheses of~\cite[Theorem 1]{BerkesPhilipp79}.
Defining 
\[
\alpha_j=16d\,T_j^{-1}\log T_j+4\lambda_j^{\frac12}T_j^d+\delta_j,
\]
as in~\cite{BerkesPhilipp79}, we have 
$\alpha_j\ll j^{-\gamma}\log j+j^{(d+\frac12)\gamma-\frac14 Q}$,
which is summable provided $1<\gamma<\frac{\frac14 Q-1}{d+\frac12}$.
We take $\gamma$ slightly larger than $1$
and $Q$ slightly larger than $4d+6$
so that $\alpha_j\ll j^{-(1+\epsilon)}$.

Applying~\cite[Theorem 1]{BerkesPhilipp79}, we conclude that 
(passing to a richer probability space)
there is a sequence of i.i.d. random
variables $Y_j$ with distribution $N(0,\Sigma)$ such that
\[
|X_j-Y_j|\ll j^{-(1+\epsilon)}\quad\text{a.e.}
\]
Let $W(t)$ be a Brownian motion with covariance $\Sigma$
and define $Y_j^*=[j^Q]^{-\frac12}(W(h_j)-W(h_{j-1}))$ where 
$h_j=\sum_{n=1}^j [n^Q]\sim j^{1+Q}$.
Then $\{Y_j\}=_d\{Y_j^*\}$ and without loss (after passing to a richer probability space), $Y_j=Y_j^*$.
We have
\begin{align*}
\SMALL
\sum_{j\le M}y_j & \SMALL=\sum_{j\le M}[j^Q]^{\frac12}X_j
=\sum_{j\le M}[j^Q]^{\frac12}(X_j-Y_j)+ \sum_{j\le M}W(h_j)-W(h_{j-1}) \\ &
\SMALL =\sum_{j\le M}[j^Q]^{\frac12}(X_j-Y_j)+ W(h_M).
\end{align*}
Now
\[
\SMALL \sum_{j\le M}[j^Q]^{\frac12}(X_j-Y_j)\ll 
\sum_{j\le M}j^{\frac12 Q}\alpha_j \ll 
\sum_{j\le M}j^{\frac12 Q-1} \ll M^{\frac12 Q}
\ll N^{\frac12 Q/(1+Q)}.
\]
If $h_M>N$, then
$h_M-N\le P_M-P_{M-1}\ll M^Q$, whereas if
$h_M<N$ then $N-h_M<P_M-h_M=\sum_{j\le M}[j^\alpha]\ll M^{1+\alpha}$.
Hence $h_M-N\ll N^{\max\{Q,1+\alpha\}/(1+Q)}$.
Taking $\alpha$ small, we obtain
$W(h_M)=W(N)+O(N^{\max\{\frac12 Q,\frac12\}/(1+Q)+\epsilon})$.
Combining these estimates with Theorem~\ref{thm-approx} we obtain
\[
\SMALL\sum_{n\le N}\eta_n=
\sum_{n\le N}\eta_n- \sum_{j\le M}y_j+ \sum_{j\le M}y_j
=W(N)+O(N^{\max\{\frac1p+\frac12 Q,\frac12\}/(1+Q)+\epsilon}).
\]
Taking $Q$ slightly larger than $4d+6$ yields the required result.
\qed

\section{ASIP for Gibbs-Markov maps}
\label{sec-GM}

In this section we prove the ASIP for weighted Lipschitz $\R^d$-valued 
observables of Gibbs-Markov maps.
Roughly speaking, these are uniformly expanding maps with countably
many inverse branches and good distortion properties, and have been
studied extensively in~\cite{Aaronson}.  We derive Theorem~\ref{thm-A}
as a consequence.

\subsection{Gibbs-Markov maps}
Let $(\Lambda,m)$ be a Lebesgue space with a countable
measurable partition $\alpha$.
Without loss, we suppose that all partition elements $a\in\alpha$ have
$m(a)>0$.
Recall that a measure-preserving transformation $F:\Lambda\to \Lambda$ is a
{\em Markov map} if $Fa$ is a union of elements  of $\alpha$
and $F|_a$ is injective for all $a\in\alpha$.
Define $\alpha'$ to be the coarsest partition of $\Lambda$ such that
$Fa$ is a union of atoms in $\alpha'$ for all $a\in\alpha$.
(So $\alpha'$ is a coarser partition than $\alpha$.)
If $a_0,\ldots,a_{n-1}\in\alpha$, we define the £$n$-cylinder
$[a_0,\ldots,a_{n-1}]=\cap_{i=0}^{n-1}F^{-i}a_i$.
It is assumed that $F$ and $\alpha$ separate points in $\Lambda$
(if $x,y\in \Lambda$ and $x\neq y$, then for $n$ large enough there
exist distinct $n$-cylinders that contain $x$ and $y$).

Let $0<\beta<1$.
We define a metric $d_\beta$ on $\Lambda$ by $d_\beta(x,y)=\beta^{s(x,y)}$
where $s(x,y)$ is the greatest integer $n\ge0$ such that
$x,y$ lie in the same $n$-cylinder.
Define $g=JF^{-1}=\frac{dm}{d(m\circ F)}$ and set
$g_k=g\,g\circ F\,\cdots\,g\circ F^{k-1}$.

A Markov map $F$ is {\em topologically mixing} if for all $a,b\in\alpha$ there
exists $N\ge1$ such that $F^na\cap b\neq\emptyset$ for all $n\ge N$.
A Markov map $F$ is {\em Gibbs-Markov} if
\begin{itemize}
\item[(i)] {\em Big images property:}   There exists $c>0$ such that
$m(Fa)\ge c$ for all $a\in\alpha$.
\item[(ii)] {\em Distortion:}  $\log g|_a$ is Lipschitz with respect to
$d_\beta$ for all $a\in\alpha'$.
\end{itemize}

Let $\alpha_0^{k-1}$ denote the partition of $\Lambda$
into length $k$ cylinders $a=[a_0,\ldots,a_{k-1}]$.
The following result of~\cite{AaronsonDenker01}
is stated explicitly in~\cite[Lemma~2.4(b)]{MN05}.

\begin{lemma} \label{lem-hyp1}
Let $F$ be a topologically mixing Gibbs-Markov map.   Then
$\bigl|m(a\cap F^{-(N+k)}b)-m(a)m(b)\bigr|\le C\tau^Nm(a)m(b)^{1/2}$
for all $a\in\alpha_0^{k-1}$ and all measurable $b$.\qed
\end{lemma}

\subsection{Weighted Lipschitz observations}

Let $p\in[1,\infty)$.  We fix a sequence of weights $R(a)>0$ satisfying
$|R|_p=(\sum_{a\in\alpha} m(a) R(a)^p)^{1/p} < \infty$.
Given $\Phi:\Lambda\to\R$ continuous, 
define $|\Phi|_\beta$ to be the Lipschitz constant of $\Phi$ with respect
to the metric $d_\beta$.  Let
$\|\Phi\|_\infty = \sup_{a\in\alpha} |\Phi1_a|_\infty/R(a)$,
$\|\Phi\|_\beta=\sup_{a\in\alpha} |\Phi1_a|_\beta/R(a)$.
Let ${\cal B}$ consist of the space of weighted Lipschitz functions
with \mbox{$\|\Phi\|=\|\Phi\|_\infty +\|\Phi\|_\beta<\infty$}.
Note in particular that $R\in {\cal B}$ and $\|R\|=1$.
We have the embeddings ${\rm Lip}\subset {\cal B}\subset L^p \subset L^1$,
where ${\rm Lip}$ is the space of (globally) Lipschitz functions.

\begin{lemma} \label{lem-hyp4}
Let $\Phi\in \mathcal{B}$ with $\int_\Lambda \Phi=0$.  
Then  $|\Phi-E(\Phi|\alpha_0^{k-1})|_p\le \|\Phi\|_\beta |R|_p\beta^k$
for all $k\ge1$.
\end{lemma}

\begin{proof} (cf.~\cite[Lemma~2.4(a)]{MN05})
Note that $E(\Phi|\alpha_0^{k-1})$ is constant on partition elements 
$a\in\alpha_0^{k-1}$ with value $\frac{1}{m(a)}\int_a\Phi\,dm$,
and that $|\Phi1_a-\frac{1}{m(a)}\int_a\Phi\,dm|_\infty\le |\Phi1_a|_\beta\diam_\beta(a)
\le \|\Phi\|_\beta R(a)\beta^k$.   Hence,
$|\Phi-E(\Phi|\alpha_0^{k-1})|_p^p\le
(\|\Phi\|_\beta\beta^k)^p\sum_{a\in\alpha_0^{k-1}}R(a)^p m(a)
= (\|\Phi\|_\beta\beta^k|R|_p)^p$.
\end{proof}

\subsection{A maximal inequality}

Given a measure-preserving transformation $F:\Lambda\to \Lambda$ of a 
probability space $(\Lambda,m)$, the 
transfer (Perron-Frobenius) operator $L:L^1\to L^1$ is defined by
$\int_\Lambda L\Phi\,\Psi\,dm=\int_\Lambda \Phi\,\Psi\circ F\,dm$
for all $\Phi\in L^1$, $\Psi\in L^\infty$.
This restricts to an operator on $L^p$, $1\le p\le\infty$.

\begin{lemma}   \label{lem-hyp3}
Let $\Phi\in L^p(\Lambda)$, $1\le p<\infty$ with $L\Phi=0$.  Then 
$\bigl|\,\max_{0\le \ell\le N-1}\bigl|\sum_{n=\ell}^{N}\Phi\circ F^n\bigr|\,\bigr|_p\le CN^{\frac12}$.
\end{lemma}

\begin{proof}
Note that $L=E(\cdot|F^{-1}\mathcal{M})$ where $\mathcal{M}$ is the underlying 
$\sigma$-algebra.  By hypothesis the sequence $\{\Phi\circ F^n;n\ge0\}$ 
is a reverse martingale difference sequence.
Passing to the natural extension we obtain an $L^p$
martingale difference sequence $\{w_n;n\in\Z\}$ such that $\Phi\circ F^n=w_{-n}$.
By Burkholder's inequality~\cite{Burkholder73}\footnote{
  This follows from~\cite[eqns (1.4) and (3.3)]{Burkholder73} and is 
  stated explicitly in~\cite[Eq.~1]{PeligradUtevWu05}.  },
we have $|\max_{1\le k\le N}|\sum_{n=0}^kw_n|\,|_p\le CN^{\frac12}$.
Setting $\ell=N-k$ and using stationarity,
\[
\SMALL\max_{0\le\ell\le N-1}|\sum_{\ell}^N\Phi\circ F^n|
=_d
\max_{0\le\ell\le N-1}|\sum_{-N+\ell}^0\Phi\circ F^n|
=
\max_{1\le k\le N}|\sum_0^kw_n|,
\]
proving the result.
\end{proof}

\subsection{Quasicompactness and the central limit theorem}

Let $F:\Lambda\to\Lambda$ be a topologically mixing Gibbs-Markov map with 
transfer operator $L:L^1\to L^1$.
It is well-known~\cite{Aaronson,MN05} that $L$ restricts
to a bounded operator on weighted Lipschitz observables $\Phi\in\mathcal{B}$
and $L(\mathcal{B})\subset {\rm Lip}$.
Moreover $L:\mathcal{B}\to\mathcal{B}$ is {\em quasicompact}:
$L1=1$ and the spectral radius of $L$ restricted to 
$\mathcal{B}_0=\{\Phi\in \mathcal{B}:\int_\Lambda \Phi\,dm=0\}$
is strictly less than $1$.

We define $\mathcal{B}^d$ to the be the space of $\R^d$ weighted
Lipschitz observables, so $\Phi=(\Phi_1,\dots,\Phi_d)\in\mathcal{B}^d$
if and only if $\Phi_i\in\mathcal{B}$ for $i=1,\dots,d$.
Similarly, we define $\mathcal{B}_0^d$.   We suppress the superscript
for spaces such as $L^p$ and ${\rm Lip}$ relying on the context.

\begin{prop}   \label{prop-ker}
Suppose that $\Phi\in\mathcal{B}_0^d$.
Then there exists $\Psi\in\mathcal{B}_0^d$ and $\chi\in L^\infty$ such that
$\Phi=\Psi+\chi\circ F-\chi$ and $L\Psi=0$.
\end{prop}

\begin{proof}   (cf.~\cite[Proof of Corollary~2.3(c)]{MN05})
Define $\chi=\sum_{j=1}^\infty L^j\Phi$.  This converges in
$\mathcal{B}_0^d$ since the spectral radius of $L$ is less than $1$.
Since $L(\mathcal{B}^d)\subset {\rm Lip}$, we have $\chi\in L^\infty$.
By construction, $L\Psi=0$.
\end{proof}

Suppose that $p\ge2$.
Let $\Phi\in\mathcal{B}_0^d\subset L^2$ and assume that $L\Phi=0$.   
Let $S_N=\sum_{n\le N}\Phi\circ F^n$ and form the $d\times d$ matrix
$S_N S_N^T$.
We define the covariance matrix
$\Sigma=\frac1N\int_\Lambda S_NS_N^T\,dm
=\int_\Lambda \Phi \Phi^T\,dm$.

\begin{lemma} \label{lem-hyp5}
There exists $\epsilon>0$ such that
\[
\SMALL\int_\Lambda\exp(i\langle u,S_N\rangle N^{-\frac12})\,dm-
\exp(-{\SMALL\frac12}\langle u,\Sigma u\rangle)=O(N^{-\frac12})
\]
uniformly for $u\in\R^d$ satisfying $|u|\le \epsilon N^\frac12$.
\end{lemma}

\begin{proof}
We follow a standard argument establishing the central limit theorem
with error term for systems with quasicompact transfer
operator (see~\cite[Theorem~4.13]{ParryPoll90} and references therein).
Let $S^{d-1}$ denote the unit sphere in $\R^d$.
Given $u\in\R^d$, write $u=tv$ where $t\ge0$ and $v\in S^{d-1}$.
Define the twisted transfer operator $L_u:\mathcal{B}\to\mathcal{B}$ by
$L_u\Psi=L(e^{i\langle u,\Phi\rangle}\Psi)$.
Recall that $1$ is an isolated eigenvalue for $L$.
For $u$ small, the spectral radius of $L_u$ is $\exp P(u)$ where
$P$ is analytic, $P(0)=0$.
Moreover~\cite[p.~66]{ParryPoll90} 
\[
P(u)=-{\SMALL\frac12}\langle u,\Sigma u\rangle -i P_3(v)t^3 +P_4(v,t)t^4,
\]
where $P_3(v)\in\R$, $P_4(v,t)\in\C$ are analytic, and 
there exists $\epsilon>0$ such that~\cite[p.~67]{ParryPoll90} 
\[
\exp(NP(uN^{-\frac12}))-\exp(-{\SMALL\frac12}\langle u,\Sigma u\rangle)
(1-iP_3(v)t^3N^{-\frac12}) =O(N^{-1})
\]
uniformly for $|u|\le \epsilon N^\frac12$.

Now $\int_\Lambda\exp(i\langle u,S_N\rangle N^{-\frac12})\,dm=
\int_\Lambda (L_{uN^{-\frac12}})^N1\,dm$, and since the leading eigenvalue of $L_u$
is isolated there exists $\gamma\in(0,1)$ such that
\[
\SMALL\int_\Lambda\exp(i\langle u,S_N\rangle N^{-\frac12})\,dm
-\exp(NP(uN^{-\frac12})) \ll \gamma^N,
\]
uniformly for $|u|\le \epsilon N^\frac12$.  This completes the proof.
\end{proof}

\subsection{Statement and proof of ASIP for Gibbs-Markov maps}

\begin{thm} \label{thm-GM}
Suppose that $F:\Lambda\to \Lambda$ is a topologically mixing Gibbs-Markov map.
Define the Banach space $\mathcal{B}^d$ corresponding to weights
$R\in L^p$ where $p>2$.
Suppose that $\Phi:\Lambda\to\R^d$ is a mean zero observable
in $\mathcal{B}^d$ with partial sums $S_N=\sum_{n=1}^N\Phi\circ F^j$.
Then the conclusion of Theorem~\ref{thm-NUH} is valid.
\end{thm}

\begin{proof}
By Proposition~\ref{prop-ker}, $S_N=\sum_{n=1}^N\Psi\circ F^j+O(1)$ a.e.
where $L\Psi=0$.
Hence without loss we may suppose from the outset that $L\Phi=0$.

Define $\eta_n=\Phi\circ F^n$ and $\xi_n=a_n$.
Then $\eta_n=\Phi(\xi_n,\xi_{n+1},\ldots)$.
We verify the hypotheses of Theorem~\ref{thm-ASIP}.

Hypothesis~\eqref{eq-mix} follows from Lemma~\ref{lem-hyp1}
and~\eqref{eq-Lp} is immediate.
The remaining hypotheses follow from Lemmas~\ref{lem-hyp3},~\ref{lem-hyp4}
and~\ref{lem-hyp5} respectively.

This completes the proof for $d\ge2$.    The improved estimate for $d=1$ follows
from Theorem~\ref{thm-1DASIP}.  (One hypothesis is different, but it
was verified in~\cite{MN05}.)
\end{proof}

\begin{pfof}{Theorem~\ref{thm-A}}
This reduces by standard techniques to a two-sided and then one-sided subshift
of finite type.  
The latter is a special case of a Gibbs-Markov map
with finite alphabet, hence $R\in L^\infty$.
Theorem~\ref{thm-A} follows from Theorem~\ref{thm-GM} with $p=\infty$.
\end{pfof}

\section{Applications to nonuniformly hyperbolic systems}
\label{sec-app}

In this section, we prove the vector-valued ASIP for large classes
of nonuniformly hyperbolic systems.   
In Subsection~\ref{sec-NUE}, we consider nonuniformly expanding systems.
In Subsection~\ref{sec-NUH}, we consider nonuniformly hyperbolic systems, 
proving Theorems~\ref{thm-NUH} and~\ref{thm-lorentz}.
Some open problems are described in Subsection~\ref{sec-open}.

\subsection{Nonuniformly expanding systems}
\label{sec-NUE}

Let $(M,d)$ be a locally compact separable bounded
metric space with Borel probability
measure $\eta$ and let $f:M\to M$ be a nonsingular transformation for
which $\eta$ is ergodic.
Let $\Lambda\subset M$ be a measurable subset with $\eta(\Lambda)>0$.
We suppose that there is an at
most countable measurable partition $\{\Lambda_j\}$ with $\eta(\Lambda_j)>0$,
and that there exist integers $R_j\ge1$, and constants $\lambda>1$;
$C>0$ and $\gamma\in(0,1)$ such that for all $j$,
\begin{itemize}
\item[(1)]  $f^{R_j}:\Lambda_j\to \Lambda$ is a (measure-theoretic) bijection.
\item[(2)]  $d(f^{R_j}x,f^{R_j}y)\ge \lambda d(x,y)$ for all $x,y\in \Lambda_j$.
\item[(3)] $d(f^kx,f^ky)\le Cd(f^{R_j}x,f^{R_j}y)$ for all $x,y\in \Lambda_j$,
$k<R_j$.
\item[(4)] $g_j=\frac{d(\eta|_{\Lambda_j}\circ (f^{R_j})^{-1})}{d\eta|_\Lambda}$ satisfies
$|\log g_j(x)-\log g_j(y)|\le Cd(x,y)^\gamma$ for almost all $x,y\in \Lambda$.
\item[(5)] $\sum_j R_j\eta(\Lambda_j)<\infty$.
\end{itemize}
A dynamical system $f$ satisfying (1)--(5) is
called {\em nonuniformly expanding}.

Define the {\em return time function} $R:\Lambda\to \Z^+$ by 
$R|_{\Lambda_j}\equiv R_j$ and the {\em induced map}
$F:\Lambda\to \Lambda$ by $Fy=f^{R(y)}(y)$.
It is well-known that there is a unique invariant
probability measure $m$ on $M$ that is equivalent to $\eta$.

\begin{thm} \label{thm-NUE}   Let $f:M\to M$ be a nonuniformly
expanding map satisfying (1)--(5) above.
Assume moreover that $R\in L^p(\Lambda)$, $p>2$.
Let $\phi:M\to\R^d$ be a mean zero H\"older observation
with partial sums $S_N=\sum_{n=1}^N\phi\circ f^j$.
Then the conclusion of Theorem~\ref{thm-NUH} is valid.
\end{thm}

\begin{proof}   This is identical to the proof of~\cite[Theorem~2.9]{MN05}
so we just sketch the main steps.   The induced map
$F:\Lambda\to \Lambda$ is a topologically mixing Gibbs-Markov map
with respect to the partition $\alpha=\{\Lambda_j\}$.
The induced observable $\Phi:\Lambda\to\R^d$ given by
$\Phi(y)=\sum_{\ell=0}^{R(y)-1}\phi(f^jy)$ is weighted Lipschitz
and satisfies the ASIP by Theorem~\ref{thm-GM}.

If $F:\Lambda\to \Lambda$ were the first return map, then the result would
follow immediately from~\cite[Theorem~4.2]{MT04} (see 
also~\cite[Theorem~B.1]{MN05}).   The general result is proved by passing
to a Young tower~\cite{Young99} which is a Markov extension of $f:M\to M$
for which $F$ is the first return map.
\end{proof}

\begin{rmk}
(a) The regularity assumption on $\phi$ in Theorem~\ref{thm-NUE}
can be replaced by the more general assumption that the induced
observable $\Phi$ is weighted Lipschitz (with respect to the 
metric defined on the Gibbs-Markov system $\Lambda$).
\\[.75ex]
(b) A similar result holds for nonuniformly expanding 
semiflows~\cite[Corollary~2.12]{MN05}.
\end{rmk}

Naturally, Theorem~\ref{thm-NUE} includes uniformly expanding and
piecewise expanding maps where the partition is finite (with $p=\infty$).
Further examples of nonuniformly expanding maps to which Theorem~\ref{thm-NUE}
applies include Alves-Viana maps,
Liverani-Saussol-Vaienti (Pomeau-Manneville maps),
multimodal maps,
and circle maps with a neutral fixed point, see~\cite[Section~4]{MN05}.

\subsection{Nonuniformly hyperbolic systems}
\label{sec-NUH}

As was the case in~\cite{MN05}, the results in this paper apply
to dynamical systems that are {\em nonuniformly hyperbolic 
in the sense of Young~\cite{Young98}} with return time
function $R\in L^p$, $p>2$.    

Let $f:M\to M$ be a diffeomorphism (possibly with singularities) defined on a
Riemannian manifold $(M,d)$.   We assume from the start that $f$ preserves
a ``nice'' probability measure $m$ (one of
the conclusions in Young~\cite{Young98} is that $m$ is a SRB measure).

Fix a subset $\Lambda\subset M$ and a family of subsets of $M$ called
``stable disks'' $\{W^s\}$ that are disjoint and cover $\Lambda$.
The stable disk containing $x$ is labelled $W^s(x)$.

\begin{itemize}
\item[(A1)]  There is a partition $\{\Lambda_j\}$ of $\Lambda$ and
integers $R_j\ge1$ such that $f^{R_j}(W^s(x))\subset W^s(f^{R_j}x)$
for all $x\in \Lambda_j$.
\end{itemize}

Define the return time function $R:\Lambda\to\Z^+$ by $R|_{\Lambda_j}=R_j$
and the induced map
$F:\Lambda \to \Lambda$ by $F(x)=f^{R(x)} (x)$.
Form the discrete suspension map $\hat f:\Delta\to\Delta$ where
$\hat f(x,\ell)=(x,\ell+1)$ for $\ell<R(x)-1$ and $\hat f(x,R(x)-1)=(Fx,0)$.
Define a separation time  $s:\Lambda\times \Lambda \to \N$ by defining
$s(x,x')$ to be the greatest integer $n\ge0$
such that $F^k x, F^k x'$ lie in the same partition element of $\Lambda$
for $k=0,\ldots,n$.
(If $x,x'$ do not lie in the same partition element, then we take
$s(x,x')=0$.)
For general points $p=(x,\ell),p'=(x',\ell')\in\Delta$,
define $s(p,q)=s(x,x')$ if $\ell=\ell'$ and $s(p,q)=0$ otherwise.
This defines a separation time $s:\Delta\times\Delta\to\N$.
The projection $\pi:\Delta\to M$, $\pi(x,\ell)=f^\ell x$,
satisfies $\pi f=\hat f\pi$.

\begin{itemize}
\item[(A2)]   There is a distinguished ``unstable leaf''
$W^u\subset\Lambda$ such that
each stable disk intersects $W^u$ in precisely one point, and
there exist constants $C\ge1$, $\alpha\in(0,1)$ such that
\begin{itemize}
\item[(i)] $d(f^nx,f^ny)\le C\alpha^n$, for all $y\in W^s(x)$, all $n\ge0$, and
\item[(ii)] $d(f^nx,f^ny)\le C\alpha^{s(x,y)}$ for all $x,y\in W^u$
and all $0\le n<R$.
\end{itemize}
\end{itemize}

Let $\bar\Lambda=\Lambda/\sim$ where $x\sim x'$ if $x\in W^s(x')$ and
define the partition $\{\bar\Lambda_j\}$ of $\bar\Lambda$.
We obtain a well-defined return time function  $R:\bar\Lambda\to\Z^+$ and
induced map $\bar F:\bar\Lambda\to\bar\Lambda$.
Let $\bar f:\bar\Delta\to\bar\Delta$ denote the 
quotient of $\hat f:\Delta\to\Delta$
where $(x,\ell)$ is identified with $(x',\ell')$ if $\ell=\ell'$
and $x'\in W^s(x)$.
Let $\bar\pi:\Delta\to\bar\Delta$ denote the natural projection.
The separation time on $\Delta$ drops down to a separation time on $\bar\Delta$.

\begin{itemize}
\item[(A3)]  The map $F:\bar\Lambda\to\bar\Lambda$ and partition $\{\bar\Lambda_j\}$ separate points in $\bar\Lambda$.
(It follows that $d_\theta(p,q)=\theta^{s(p,q)}$ defines
a metric on $\bar\Delta$ for each $\theta\in(0,1)$.)
\item[(A4)]  There exist
invariant probability measures $\hat m$ on $\Delta$ and $\bar m$ on
$\bar\Delta$ such that
\begin{itemize}
\item[(i)]  $\pi:\Delta\to M$ and $\bar\pi:\Delta\to \bar\Delta$ are
measure-preserving; and 
\item[(ii)]  $\bar f:\bar\Delta\to\bar\Delta$ is nonuniformly expanding in the 
sense of Subsection~\ref{sec-NUE} with induced map 
$\bar F:\bar\Lambda\to\Lambda$.  (Conditions (2) and (3) are automatic.)
\end{itemize}
\end{itemize}

\begin{pfof}{Theorem~\ref{thm-NUH}}
This reduces, as in the proof of~\cite[Theorem~3.4]{MN05}, 
to the ASIP for the nonuniformly expanding map $\bar f:\bar\Delta\to\bar\Delta$
and hence follows from Theorem~\ref{thm-NUE}.
\end{pfof}

\begin{rmk} Again, the regularity assumption on $\phi$ can be relaxed, and
the result extends to nonuniformly hyperbolic flows.
\end{rmk}

Large classes of billiard maps and Lorentz flows, surveyed 
in~\cite{ChernovYoung00} satisfy the vector-valued ASIP.   These include
dispersing billiards (with finite or infinite horizons) 
and the corresponding Lorentz flows (assuming finite horizons).

\begin{pfof}{Theorem~\ref{thm-lorentz}}
By periodicity, we can consider the quotient flow on the compact manifold
$M=(\T^2-\Omega)\times S^1$.   The Poincar\'e map $f:X\to X$
on the compact cross-section
$X=\partial\Omega\times(\frac{-\pi}{2},\frac{\pi}{2})$ is 
called the billiard map or collision map.
Benedicks \& Young~\cite{BenedicksYoung00} showed that $f$ is nonuniformly 
hyperbolic in the sense of Young with $R\in L^p$ for all $p>2$.
By Theorem~\ref{thm-NUH}, the vector-valued ASIP
holds for $f$ with $p=\infty$.
The collision time is uniformly bounded and piecewise H\"older, so
it follows from~\cite{MT04} that the vector-valued ASIP
holds for the Lorentz flow on $M$.  Now take as an $\R^2$-valued observable 
the velocity coordinate $v:M\to S^1$.   This is piecewise H\"older, and the 
lifted position in $\R^2$ is given by $q(T)=\int_0^T v\circ f_t\,dt$.  Finally, 
nonsingularity of the covariance matrix was proved in~\cite{BunimSinai81}.
\end{pfof}

Chernov \& Zhang~\cite{ChernovZhang05} consider
three classes of billiards with slow mixing rates.
The first and third classes are not covered by
our results since it is shown only that $R\in L^{2-\epsilon}$.   
The second class of {\em Bunimovich-type}
billiards treated in~\cite{ChernovZhang05} satisfies 
$R\in L^{3-\epsilon}$.  The vector-valued ASIP for such billiards (and the 
corresponding flows) is hence a consequence of Theorem~\ref{thm-NUH}
with error $\beta=\frac{6d+10}{12d+21}$.

As in~\cite{MN05}, Theorem~\ref{thm-NUH} also applies to 
Lozi maps and certain piecewise hyperbolic
maps, H\'enon-like maps and partially hyperbolic
diffeomorphisms with mostly contracting direction.
 
A further important class of dynamical systems is {\em singular
hyperbolic flows}~\cite{MoralesPacificoPujals99}.   Theorem~\ref{thm-NUH}
does not apply directly to such systems, but it establishes
the vector-valued ASIP (with $p=\infty$) when combined with the techniques
in Holland \& Melbourne~\cite{HMsub}.

\subsection{Open problems}
\label{sec-open}

Given a (non)uniformly hyperbolic flow $f_t$, the time-one map $f_1$ is only {\em partially hyperbolic}.   For such maps the Gibbs-Markov/suspension formalism
breaks down so the results in~\cite{MN05} and in this paper do not apply.
By different methods,
Melbourne \& T\"or\"ok~\cite{MT02} proved that the scalar ASIP
is typically valid for the time-one map of an Axiom~A flow.
They used rapid mixing properties to reduce to a reverse martingale
difference sequence.  
Following~\cite{ConzeBorgne01,FMT03}, this leads to the ASIP
in reverse time and hence forwards time (since the class of such flows
is closed under time reversal).
Similarly, the scalar ASIP for the time-one map of the planar periodic
Lorentz gas with finite horizons is typically valid
(since the flow is typically rapid mixing~\cite{Mapp} and the class
of flows is closed under time reversal).

\begin{problem}  \label{prob-1}
Prove that the vector-valued ASIP holds (at least typically) for
time-one maps of Axiom~A flows and/or planar periodic
Lorentz gas with finite horizons.
\end{problem}

Generally speaking, the hypotheses for a nonuniformly hyperbolic system are
not time-symmetric so~\cite{ConzeBorgne01,FMT03,MT02} fails.   

\begin{problem}  \label{prob-2}
Obtain results on the scalar ASIP for time-one maps 
of nonuniformly hyperbolic flows and/or singular hyperbolic flows.
\end{problem}

\begin{rmk}
(a) The Banach space-valued ASIP of~\cite{Berger90} applies to 
Problem~\ref{prob-1},
with the caveats mentioned in Remark~\ref{rmk-Berger}(a).
In particular, the $d$-dimensional functional LIL is typically valid.
These results do not apply to Problem~\ref{prob-2}.
\\[.75ex]
(b) The (vector-valued) functional central limit theorem is typically valid
in Problems~\ref{prob-1} and~\ref{prob-2} (combining the arguments 
in~\cite[Section~3.3]{FMT03} and~\cite{MT02}).
\end{rmk}

\appendix

\section{Scalar ASIP with error term}
\label{sec-1}

In this appendix, we prove a scalar ASIP using 
martingale approximation and the method of Strassen~\cite{Strassen67}.
This is precisely
the result~\cite[Theorem~7.1]{PhilippStout75} used in~\cite{MN05},
but our purpose here is to obtain a better error term under assumptions 
appropriate for dynamical systems.
This improves Theorem~\ref{thm-ASIP} when $d=1$.

We assume the conditions of Section~\ref{sec-ASIP} 
except that~\eqref{eq-CLT} is replaced by
\begin{align} \label{eq-variance}
\SMALL E(\sum_{n\le N}\eta_n)^2=N+O(N^{1/2}).
\end{align}
Define $\{y_j\}$ as in Section~\ref{sec-block}.   
In particular, Theorem~\ref{thm-approx} is unchanged.

\subsection*{Law of large numbers for $y_j^2$}


\begin{lemma} \label{lem-LLNE}
Let 
$\gamma=\max\bigl\{\frac{\frac12 Q}{1+Q}, \frac{\frac12+\alpha}{1+Q} \bigr\}$.
Then
$\sum_{j\le M_N}Ey_j^2 = N +O(N^{\frac12+\gamma})$.
\end{lemma}

\begin{proof}
(cf.~\cite[Lemma~7.3.5]{PhilippStout75})
By~\eqref{eq-variance}, 
 $E(\sum_{n\le N}\eta_n)^2=a_N^2$ where $a_N^2=N(1+O(N^{-1/2}))$.
Write 
\[
\sum_{n\le N}\eta_n= \sum_{n\le P_M}\eta_n - \sum_{n=N+1}^{P_M}\eta_n
= \sum_{n\le P_M}(\eta_n-\eta_{\ell,n})+ \sum_{j\le M}y_j+\sum_{j\le M}z_j
- \sum_{n=N+1}^{P_M}\eta_n.
\]
Then
\begin{align*}
\SMALL\bigl\|\sum_{j\le M}y_j\bigr\|_2-a_N &=
\SMALL\bigl\|\sum_{j\le M}y_j\bigr\|_2-\bigl\|\sum_{n\le N}\eta_n\bigr\|_2 \\
& \SMALL\le \bigl\|\sum_{j\le M}z_j\bigr\|_2+\bigl\|\sum_{n\le P_M}(\eta_n-\eta_{\ell n})\bigr\|_2+\bigl\|\sum_{n=N+1}^{P_M}\eta_n\bigr\|_2.
\end{align*}
By Proposition~\ref{prop-local}, 
$\|\sum_{n\le N}(\eta_n-\eta_{\ell n})\|_2\ll1$.
By the proof of Lemma~\ref{lem-blocks2},
$\|\sum_{j\le M}z_j\|_2^2\ll M^{1+2\alpha}$ and so
$\|\sum_{j\le M}z_j\|_2\ll N^{(\frac12+\alpha)/(1+Q)}$.
By stationarity and~\eqref{eq-variance}, 
$\|\sum_{n=N+1}^{P_M}\eta_n\|_2^2 =
\|\sum_{n\le P_M-N}\eta_n\|_2^2\ll P_M-N\ll N^{Q/(1+Q)}$.
Hence $\|\sum_{j\le M}y_j\|_2=a_N+O(N^{\gamma})$ and 
$E(\sum_{j\le M}y_j)^2=N+O(N^{\frac12+\gamma})$.
Also, as in the proof of Lemma~\ref{lem-blocks2}, 
$\sum_{i\neq j} E y_iy_j \ll 1$.~
\end{proof}

\begin{cor} \label{cor-LLNE}
Let $\beta=\max\bigl\{\frac{\frac14+\frac12 Q}{1+Q}, 
\frac{\frac12+\frac12\alpha+\frac14 Q}{1+Q} \bigr\}$.  Then
$\sum_{j\le M_N}Ey_j^2 = N +O(N^{2\beta})$.~\qed
\end{cor}


\begin{lemma} \label{lem-LLN}
Let $\beta=(\frac34-\frac{p}{8}+\frac12 Q)/(1+Q)$ for $2<p\le4$
(and $\beta=(\frac14+\frac12 Q)/(1+Q)$ for $p>4$).
Then for any $\epsilon>0$,
$\sum_{j\le M_N}y_j^2-Ey_j^2\ll N^{2\beta+\epsilon} \enspace\text{a.e.}$
\end{lemma}

\begin{proof}
The value of $\epsilon$ below  may change from line to line.
Define
\[
w_j=\begin{cases} y_j^2-Ey_j^2, & |y_j^2-Ey_j^2|\le j^{1+Q+\epsilon} \\
0, & \text{otherwise} \end{cases}
\]
Note that $P(w_j\neq y_j^2-Ey_j^2)=P(|y_j^2-Ey_j^2|> j^{1+Q+\epsilon})
\le 2\|y_j\|_2^2/j^{1+Q+\epsilon}\ll j^{-(1+\epsilon)}$
which is summable, so by Borel-Cantelli $w_j$ fails to coincide with
$y_j^2-Ey_j^2$ only finitely often.  Hence it suffices to estimate
$\sum_{j\le M}w_j$.    We do this by estimating $\sum_{j\le M}\tilde w_j$ and
$\sum_{j\le M}Ew_j$ where $\tilde w_j=w_j-Ew_j$.

Again $\sum_{i\neq j}E\tilde w_i\tilde w_j\ll 1$.  Also, 
$E\tilde w_j^2\le |\tilde w_j^{p/2}|_1\|\tilde w_j^{2-p/2}\|_\infty\ll 
\|y_j\|_p^p\|w_j\|_\infty^{2-p/2} \ll j^{R-1}$, where 
$R=3-\frac{p}{2}+2Q+\epsilon$.
Hence $E(\sum_{j=m}^n \tilde w_j)^2 \ll n^R-m^R$.
By Lemma~\ref{lem-GK}, $\sum_{j\le M} \tilde w_j\ll 
M^{\frac12 R}\le N^{\frac12 R/(1+Q)}$ for any $\epsilon>0$.

Let $A=\{|y_j^2-Ey_j^2|>j^{1+Q+\epsilon}\}$.  Then
\begin{align*}
Ew_j & =-E\Bigl\{(y_j^2-Ey_j^2)I_A\Bigr\}\ll \|y_j^2-Ey_j^2\|_{p/2} 
\|1_A\|_{p/(p-2)} \\
& \ll \|y_j\|_p^2 \|1_A\|_{p/(p-2)}
\ll j^{Q-(1+\epsilon))(p-2)/p}=j^{S-1}.
\end{align*}
where $S=\frac2p-\epsilon+Q$.
Hence $\sum_{j\le M} Ew_j\ll M^S$.
Thus, it suffices that 
$2\beta=\max\{\frac12 R/(1+Q),S/(1+Q)\} =\frac12 R/(1+Q)$.
\end{proof}

\subsection*{Martingale approximation}

\begin{lemma} \label{lem-mart}
Set $\mathcal{L}_j=\mathcal{L}_1^j=\sigma\{y_1,\dots,y_j\}$.
There is a martingale difference sequence
$\{Y_j,\mathcal{L}_j\}$ such that $y_j=Y_j+u_j-u_{j+1}$,
where $\|u_j\|_q\ll \tilde\tau^{j^\alpha}$
for all $2<q<p$.
\end{lemma}

\begin{proof}  (cf.~\cite[Lemma~7.4.1]{PhilippStout75})
Define $u_j=\sum_{k=0}^\infty E(y_{j+k}|\mathcal{L}_{j-1})$.
We estimate $\|E(y_{j+k}|\mathcal{L}_{j-1})\|_q$ which we write for convenience
as $\|E(y|\mathcal{L})\|_q$.
Note that
\begin{align*}
E|E(y|\mathcal{L})|^q & =E\Bigl\{E(y|\mathcal{L})E(y|\mathcal{L})|E(y|\mathcal{L})|^{q-2}\Bigr\}
=E\Bigl\{E\bigl\{yE(y|\mathcal{L})|E(y|\mathcal{L})|^{q-2}|\mathcal{L}\bigr\}\Bigr\} \\ 
&= E\bigl\{yE(y|\mathcal{L})E(y|\mathcal{L})^{q-2}\bigr\}.
\end{align*}
Write $\frac1q+\frac1s=1$.    Then $\frac1p+\frac1s<1$, so
by Lemma~\ref{lem-Dav},
\[
E|E(y|\mathcal{L})|^q \le \|y\|_p\|E(y|\mathcal{L})^{q-1}\|_s\tilde\tau^{(j+k)^\alpha}.
\]
Note that $\|E(y|\mathcal{L})^{q-1}\|_s=(E|E(y|\mathcal{L})|^q)^{1-\frac1q}$, and so
dividing both sides by this yields
$\|E(y|\mathcal{L})\|_q\le \|y\|_p\tilde\tau^{(j+k)^\alpha}$.
Since $\sum_{k=0}^\infty\tilde\tau^{(j+k)^\alpha}
\ll j\tilde\tau^{j^\alpha}$, the estimate for $\|u_j\|_q$ follows
(increasing $\tilde\tau$ slightly).

At the same time, it follows immediately that
$\sum_{k=0}^\infty |E(y_{j+k}|\mathcal{L}_j)|_1<\infty$ which guarantees
(see eg.~\cite[Lemma~2.1]{PhilippStout75}) that $Y_j$ is a martingale
difference sequence.
\end{proof}

\begin{cor} \label{cor-mart}
$\sum_{j\le M_N}(y_j-Y_j) \ll 1\enspace\text{a.e.}$
and $\sum_{j\le M_N}(y_j^2-Y_j^2) \ll N^{\frac12}\enspace\text{a.e.}$
\end{cor}

\begin{proof}
We have $\sum_{j\le M}(y_j-Y_j)=u_1-u_{M+1}$,
and hence certainly
$|\sum_{j\le M}(y_j-Y_j)|\le \sum_{j\ge1}|u_j|$.
By Lemma~\ref{lem-mart}, 
$\sum_{j\ge1}|u_j|_1<\infty$ so that $\sum_{j\ge1}|u_j|<\infty$
a.e. proving the first statement.

Set $v_j=u_j-u_{j+1}$.
Then $Y_j^2-y_j^2=v_j^2-2y_jv_j$.
Now $E\sum_{j=1}^\infty v_j^2
=\sum_{j=1}^\infty Ev_j^2\ll\sum_{j=1}^\infty\tilde\tau^{2j^\alpha}<\infty$ by 
Lemma~\ref{lem-mart}.
Hence $\sum_{j\le M}v_j^2\le\sum_{j=1}^\infty v_j^2<\infty$ a.e.
Finally, $\sum_{j\le M}y_jv_j\le (\sum_{j\le M}y_j^2)^{\frac12}(\sum_{j\le M}v_j^2)^\frac12\ll   (\sum_{j\le M}y_j^2)^{\frac12}\ll N^\frac12$ by 
Theorem~\ref{thm-approx}.
\end{proof}

\begin{lemma} \label{lem-Chow}
Let $X_k$ be a martingale difference sequence
and $m\in(1,2]$.   Suppose that $b_1<b_2<\cdots\to\infty$.  If
$\sum_{k\le n}b_k^{-m}E|X_k|^m<\infty$, then $\sum_{k\le n}X_k=o(b_n)$ a.e.
\end{lemma}

\begin{proof}
For $m=2$, this is explicit in~\cite[p.~238]{Feller66}.
For $m\in(1,2)$ it follows from a standard martingale result,
Chow~\cite{Chow65}, combined with Kronecker's lemma; this is 
implicit in the proof of~\cite[Lemma~7.4.4]{PhilippStout75}.
\end{proof}

\begin{lemma} \label{lem-mart2}
Let $\beta=(\frac1p+\frac12 Q)/(1+Q)$ for $2<p\le4$
(and $\beta=(\frac14+\frac12 Q)/(1+Q)$ for $p>4$).
Then for any $\epsilon>0$,
$\sum_{j\le M_N}(E(Y_j^2|\mathcal{L}_{j-1})-Y_j^2)\ll N^{2\beta+\epsilon}
\enspace\text{a.e.}$
\end{lemma}

\begin{proof}   
Define $R_j=E(Y_j^2|\mathcal{L}_{j-1})-Y_j^2$.
Suppose that $\gamma>2/p+Q$ and choose $q<p$ so that
$\gamma>2/q+Q$.
Then
\[
(j^\gamma)^{-q/2}E|R_j|^{q/2}\ll  j^{-\gamma q/2}E|Y_j|^q \ll 
j^{-(\gamma-Q)q/2},
\]
hence $\sum_{j=1}^\infty (j^\gamma)^{-q/2}E|R_j|^{q/2}<\infty$.
Note that $\frac{q}{2}\in(1,2]$ and $R_j$ is a martingale difference sequence,
so it follows from Lemma~\ref{lem-Chow} that
$\sum_{j\le M}R_j\ll M^{\gamma}$ and the result follows from 
Proposition~\ref{prop-M}.
\end{proof}

We now apply Strassen's martingale version of the Skorokhod
embedding~\cite{Strassen67}.   There exist non-negative random variables
$T_j$ such that the sequences
$\bigl\{\sum_{j\le M}Y_j,\,M\ge1\bigr\}$ and
$\bigl\{W(\sum_{j\le M}T_j),\,M\ge1\bigr\}$
are equal in distribution.   

\begin{prop} \label{prop-T}
For $\beta$ as in Corollary~\ref{cor-LLNE}, Lemma~\ref{lem-LLN},
Lemma~\ref{lem-mart2} and any $\epsilon>0$,
$\sum_{j\le M_N} T_j -N \ll N^{2\beta+\epsilon}\enspace\text{a.e.}$
\end{prop}

\begin{proof}
Let $\mathcal{A}_M=\sigma\{W(t):0\le t\le \sum_{j\le M}T_j\}$,
so $\mathcal{L}_M\subset\mathcal{A_M}$.  
Each $T_j$ is $\mathcal{A}_j$-measurable,
$E(T_j|\mathcal{A}_{j-1})=E(Y_j^2|\mathcal{L}_{j-1})$ a.e., and
$E T_j^p\ll E|Y_j|^{2p}$.
In particular, the argument in Lemma~\ref{lem-mart2}
implies that 
\begin{align} \label{eq-T}
\SMALL\sum_{j\le M}(T_j-E(T_j|\mathcal{A}_{j-1}))\ll N^{2\beta}
\enspace\text{a.e.}
\end{align}
Now write
\[
\SMALL\sum T_j-N=\sum(T_j-E(T_j|A_{j-1}))+\sum(E(Y_j^2|\mathcal{L}_{j-1})-Y_j^2)
+\sum Y_j^2-N.
\]
The result follows from~\eqref{eq-T},
Corollaries~\ref{cor-LLNE} and~\ref{cor-mart}, 
and Lemmas~\ref{lem-LLN} and~\ref{lem-mart2}.
\end{proof}

\begin{thm} \label{thm-1DASIP}
Let $\beta=\frac{1}{2p}+\frac14$ for $2<p\le4$
and $\beta=\frac38$ for $p>4$.
For any $\epsilon>0$,
$\sum_{n\le N}\eta_n = W(N)+O(N^{\beta+\epsilon})\enspace\text{a.e.}$
\end{thm}

\begin{proof}
By Theorem~\ref{thm-approx} and Corollary~\ref{cor-mart}, 
it suffices to prove that
$\sum_{j\le M}Y_j=W(N)+O(N^{\beta+2\epsilon})$.
Equivalently,
$W(\sum_{j\le M}T_j)=W(N)+O(N^{\beta+2\epsilon})$.
Hence the result follows from Proposition~\ref{prop-T}.
\end{proof}

\paragraph{Acknowledgments}  
The research of was supported in part by EPSRC Grant EP/D055520/1 and
a Leverhulme Research Fellowship (IM), and by NSF grant DMS-0244529 (MN).
IM is grateful to the University of Houston for hospitality during
part of this project, and for the use of
e-mail given that pine is inadequately supported on the University
of Surrey network.
MN would like to thank the University of Surrey for
hospitality during part of this research.

\end{document}